# Partially specified prior

K. Govindaraju and G. Jones

*Institute of Fundamental Sciences*

*Massey University*

*New Zealand*

*e-mail:* `k.govindaraju@massey.ac.nz`; `g.jones@massey.ac.nz`

**Abstract:** This note introduces the concept of a partially specified prior distribution for certain *post hoc* inference problems, where a finite population is sampled once in order to make a decision on the presence or complete absence of some attribute. If the decision is made to accept complete absence, a probability statement may be required that the population is indeed free of the attribute. A partially specified prior is shown to be advantageous in making such statements realistic and useful.

## 1. Introduction

Consider a finite population of size $N$ in which $A$ of the elements posses a certain characteristic or attribute. Often the desired value for $A$ is zero, and the decision problem is $\{A = 0\}$ versus $\{A > 0\}$. For example, certification may be required that a population of animals is free from an infectious disease, as against having one or more infected animals. Alternatively, in a quality inspection setting, we may require $\{A \leq k\}$ where $k \ll N$.

Suppose that a simple random sample of size $n$ without replacement is drawn, from which $a$ of the sampled elements are found to possess the attribute. When $a \geq 1$, the inference on whether $A = 0$ or not is obvious. When $a = 0$, a decision may be made to accept the hypothesis that $A = 0$, or $A \leq k$. The consequent issue after making this decision is to ascertain the probability that $A = 0$ or $A \leq k$. One indirect way of answering this question is to construct confidence or credible limits for $A$; see, for example, see [7] and [3]. The problem of constructing a uniformly most accurate $100(1-\alpha)\%$ upper confidence bound for $A$ is addressed in [9]. It was also remarked in [9] that if one is to be 99% "confident" that the population is free from the set attribute, then the sampling fraction $n/N$ must be at least 99%. The use of Bayesian methods was suggested by [2]. It was also commented that $\Pr\{A = 0|a = 0\}$ is of practical interest.

The rest of this note is organised as follows. In Section 2, we review the results given in [2] and discuss the determination of $\Pr\{A = 0|a = 0\}$. In section 3, the concept of partially specified prior is introduced and illustrated. The final section provides the conclusions.

## 2. Determination of $\Pr\{A = 0|a = 0\}$

For a given discrete prior distribution $p(A)$ on the support $Z_{N+1} = \{0, 1, 2..., N\}$, the following solution is provided in [2]):

$$\Pr\{A = 0 | a = 0\} = \frac{\binom{N}{n} p(0)}{\sum_{A=0}^{N-n} \binom{N-A}{n} p(A)} \tag{2.1}$$

In the case of a uniform prior $p(A) = 1/(N+1)$, see [2],

$$\Pr\{A = 0 | a = 0\} = (n+1)/(N+1) \tag{2.2}$$

Equation (2.2) suggests that the sampling fraction $n/N$ must be very large to achieve a large value for $\Pr\{A = 0 | a = 0\}$ such as 0.95 or 0.99. This is a consequence of the employment of a "non-informative" prior distribution.

A Bayesian approach for the determination of credible intervals for the fraction nonconforming, using a uniform prior with an assumed upper bound $M$ was discussed in [4]. Using their approach we obtain

$$\Pr\{A = 0 | a = 0\} = \frac{\binom{N}{n}}{\binom{N+1}{n+1} - \binom{N-M}{n+1}}$$

or equivalently

$$\Pr\{A = 0 | a = 0\} = \frac{n+1}{N+1} \left( 1 - \frac{\binom{N-M}{n+1}}{\binom{N+1}{n+1}} \right)^{-1}$$

Again the sampling fraction $n/N$ must be very large to achieve a large posterior probability, which leads to the consideration of other, perhaps more realistic, priors.

The mixed binomial prior distribution was originally proposed by [1] in a quality control context. This prior was later generalised in [8]. The basic mixed binomial prior distribution given in [1] consists of two binomials which are mixed with known probabilities of mixing. This prior, originally suggested for lot by lot sampling inspection, is of the form

$(\alpha_i, p_i)$, $(i = 1, 2)$ with $\alpha_1 + \alpha_2 = 1$

where $\alpha_i$ is the probability with which a lot has a proportion $p_i$ non-conforming. The basic idea behind the mixed binomial prior is that lots are either of acceptable quality $p_1$, under common cause variations in the production process, or of unacceptable quality $p_2\,(>p_1)$ in the presence of special causes. One advantage of the mixed binomial prior is that it takes the structure of the inspection problem into account. The main purpose of sampling inspection is to dispose the lot under consideration either as acceptable or unacceptable. The mixed prior distribution represents the two subgroups, "good" and "bad", of the superpopulation (production process) from which the lot under consideration is formed.

A practical problem with the mixed binomial prior is that it is difficult to specify in full. Whilst some knowledge of the characteristics of common cause variation might be expected, special causes by their very nature are unpredictable and unquantifiable. In the epidemiological context, some information may be forthcoming about background rates of infection, but not about sporadic outbreaks of disease. There may in fact be many kinds of "bad" lot, requiring a prior with many components. We note however that full specification of the prior is not necessary if one considers the two-stage nature of the inferential procedure. First a decision is made, based on an explicit acceptance rule, as to whether the lot at hand is "good" or "bad". Subsequently, a probability statement on the number of non-conformances is only required when a lot is accepted as "good". At this point we can dispense with those components of the prior that reference "bad" lots, and work only with the component that references a "good" lot. Thus the prior only needs to be specified for the "good" lots.

This differs from the usual Bayesian inference in that, instead of updating the prior based on the data to give p(A|Data), we are using the data to make a decision and then updating the prior based on the decision, p(A|Decision). We define the term "partially specified prior" or PSP as the required component of the full prior specified for *post hoc* inference after a decision on the population of interest is made. We do not advocate that such partial specification is possible or advisable for all inference problems. This approach is suitable only for problems whose prior distributions can be thought to be mixtures of two or more individual distributions, and where an explicit decision rule has been agreed *a priori*.

Strictly speaking, the term "prior" refers to the specification of the full distribution of the parameters before obtaining data. We suggest that the use of a PSP can be interpreted as similar to the adjusting the belief structures discussed in [6]. The lack of general methods for dealing with partial specification of probability structures is also highlighted in [6].

## 3. Partially specified prior for $\Pr\{A = 0 | a = 0\}$

For decision problems on a finite population, a decision is first made based on the sample data. For example, a batch of items may be accepted if the sampled items are all conforming. The sample size for such acceptance sampling problems is not determined purely based on statistical grounds (Bayesian or otherwise). Cost and other industry standards or mutual contracts between the producer and the consumer dictate such sample sizes. If the lot is accepted after the observance of $a = 0$, then evidently $0 < \Pr\{A = 0\} < 1$ for $n < N$. The relevant post hoc inference question is what is $\Pr\{A = 0\}$ or perhaps $\Pr\{A \leq k\}$ where $k$ is a small number. To answer such questions, consideration of an appropriate prior can be related to the structure of the problem, namely that the hypothesis $A = 0$ will be accepted under the given sampling plan if $a = 0$ is observed. Since a *post hoc* probability statement is only required following acceptance of a lot, we can ignore those components of the prior related to bad quality lots, and work only with a partial specification of the prior information. We assume that some prior knowledge of belief concerning the good lots is available, perhaps in the form of a Binomial($N$, $\delta$) comprising that part of the mixed binomial prior pertaining to a good lot given in [1]. Note that the structure of the problem in hand is important to adopt a mixture prior and then separate out the relevant part of it to form the partially specified prior for post hoc inference.

Even though we focus largely on the case $a = 0$ in this note, the suggested binomial partially specified prior is a conjugate prior for all values of $a$ since

$$\begin{aligned}
P(A|a) &\propto P(a|A)p_{PS}(A) \\
&\propto \frac{A!(N-A)!}{(A-a)!(N-A-n+a)!} \frac{N!}{A!(N-A)!} \delta^A (1-\delta)^{N-A} \\
&\propto \binom{N-n}{A-a} \delta^{A-a}(1-\delta)^{N-A-n+a}
\end{aligned} \quad (3.1)$$

for $A = a, a+1, \ldots, N$, implying that the posterior distribution of $A - a$ is Binomial$(N - n, \delta)$.

The parameter $\delta$ for the partially specified prior $p_{PS}(A)$ can be elicited as the expected proportion of elements having the attribute in a good lot. Note that:

$$\begin{aligned}
P(a = 0) &= \sum_{A=0}^{N-n} P(a = 0) \times p_{PS}(a) \\
&= \sum_{A=0}^{N-n} \left( \binom{N-A}{n} \bigg/ \binom{N}{n} \right) \times \binom{N}{A} \delta^A(1-\delta)^{N-a} \\
&= \sum_{A=0}^{N-n} \binom{N-A}{n} \delta^A(1-\delta)^{N-A} \\
&= (1-\delta)^n
\end{aligned} \quad (3.2)$$

Clearly $(1 - \delta)^n$ is the probability that all $n$ items are free from the attribute interest in which $\delta$ is the probability of an element having the attribute in a very large population.

From equation (3.1), we obtain

$$\Pr\{A = 0 | a = 0\} = (1 - \delta)^{N-n} \quad (3.3)$$

As an example, consider the international Standard ISO 2859-1:1999 [5] normal inspection plan for the lot size $N = 3200$ and Acceptance Quality Limit ($AQL$) of 0.10%. The sampling plan $n = 125$ and $Ac = 0$ is listed in the ISO Standard corresponding to the desired $AQL$ and lot size. After acceptance of the lot, the belief of the consumer improves. The term $AQL$ is the *maximum* percent nonconforming that is considered satisfactory as a process average $p$. After acceptance of the lot, the improved consumer's belief is that $p \leq AQL$. Hence a plausible value for $\delta$ is $0.5AQL$. Based on this value for $\delta$, we find from equation (3.3) that $\Pr\{A = 0|a = 0\} = 0.215$ under the plan $n = 125$ and $Ac = 0$. his figure is much larger than 0.0394 found using the uniform prior (equation (2.2)).

Finally we examine briefly the case of a non-zero acceptance number, where a lot is accepted on observance of $a \leq Ac$ where $Ac > 0$. Equation (3.1) can be used to calculate the probability that the lot is of acceptable quality after acceptance of the lot namely:

$$P(A \leq k|a) = \sum_{A=0}^{k-a} \binom{N-n}{A-a} \delta^{A-a}(1-\delta)^{N-A-n+a} \qquad (3.4)$$

where $k = int\,(NAQL)$. For example, consider the normal inspection plan $n = 200$ and $Ac = 2$ corresponding to $N = 10,000$ and AQL $= 0.40\%$ from the ISO 2859-1:1999 Standard. For given $a = 1$, $k = 10000 \times 0.004 = 40$, and $\delta = 0.002$, equation (3.4) yields 99.9%.

## 4. Concluding Remarks

This paper introduces the concept of partially specified prior for *post hoc* statistical inference issues such as making probability statements about a finite population of interest after certain decisions are made concerning the population. The suggested partial specification of the prior is useful for quality control and other applications.